\documentclass[12pt]{article}
\usepackage{fullpage}
\usepackage{graphicx}
\usepackage{amssymb}
\usepackage{epstopdf}

\usepackage{fullpage}

\def\qed{\hfill
\ifhmode\unskip\nobreak\fi\quad\ifmmode\Box\else$\Box$\fi\\ }
\newtheorem{theorem}{Theorem}
\newtheorem{cor}[theorem]{Corollary}
\newtheorem{lemma}[theorem]{Lemma}
\newtheorem{defn}[theorem]{Definition}

\newtheorem{claim}[theorem]{Claim}
\newtheorem{fact}[theorem]{Fact}
\newtheorem{conj}[theorem]{Conjecture}

\title{Ore's Conjecture on color-critical graphs is almost  true}
\author{Alexandr  Kostochka\thanks{
University of Illinois at Urbana--Champaign, Urbana, IL 61801, USA and
 Sobolev Institute of Mathematics, Novosibirsk 630090, Russia. Email:
 kostochk@math.uiuc.edu.
Research of this author
is supported in part by NSF grant DMS-0965587 and by
grants 12-01-00448 and 12-01-00631 of the Russian Foundation for Basic Research. }
\and
Matthew Yancey\thanks{Department of Mathematics, University of Illinois, Urbana,
IL 61801, USA. E-mail: yancey1@illinois.edu.
Research of this author is partially supported by the
Arnold O. Beckman Research Award of the University of Illinois
at Urbana-Champaign and from National Science Foundation grant DMS 08-38434 ``EMSW21-MCTP: Research
Experience for Graduate Students.''}}

\begin{document}

\maketitle

\begin{abstract}
A graph $G$ is  $k$-{\em critical} if it has chromatic number $k$, but every
proper subgraph of $G$ is $(k-1)$--colorable.
Let $f_k(n)$ denote the minimum number of edges in an $n$-vertex $k$-critical graph.
We give a lower bound, $f_k(n)  \geq F(k,n)$, that is sharp for every $n=1\,({\rm mod }\, k-1)$.
It is also sharp for $k=4$ and every $n\geq 6$.
The result improves the classical bounds by Gallai and Dirac and subsequent bounds by Krivelevich and
Kostochka and Stiebitz.
It establishes the asymptotics of $f_k(n)$ for every fixed $k$.
It also proves that  the conjecture by Ore from 1967 that
for every $k\geq 4$ and $n\geq k+2$, $f_k(n+k-1)=f(n)+\frac{k-1}{2}(k - \frac{2}{k-1})$
holds for each $k\geq 4$ for all but at most $k^3/12$  values of $n$.
We give a polynomial-time algorithm for $(k-1)$-coloring a graph $G$ 
that satisfies $|E(G[W])| < F_k(|W|)$ for all $W \subseteq V(G)$, $|W| \geq k$.
We also present some applications of the result.
\end{abstract}

\section{Introduction}

A {\em proper $k$-coloring}, or simply $k$-{\em coloring},
 of a graph $G = (V, E)$ is a function $f:V \rightarrow \{1,2,\dots,k\}$ such that for each $uv \in E$, $f(u) \neq f(v)$.
A graph $G$ is $k$-{\em colorable} if there exists a $k$-coloring of $G$. The {\em chromatic number}, $\chi(G)$,
of a graph $G$ is the smallest $k$ such that $G$ is $k$-{colorable}. A graph $G$ is $k$-{\em chromatic} if $\chi(G)=k$.

A graph $G$ is $k$-{\em critical} if $G$ is not $(k-1)$-colorable, but  every proper subgraph of $G$ is $(k-1)$-colorable.
  Then every $k$-critical
graph has chromatic number $k$ and every $k$-chromatic graph contains
a $k$-critical subgraph.  The importance of the notion of criticality
is that problems for $k$-chromatic graphs may often be reduced to
problems for $k$-critical graphs, whose structure is more restricted. For example, every $k$-critical graph is
$2$-connected and $(k-1)$-edge-connected.

Critical graphs were first defined and used by Dirac~{\cite{D0,D02,D03}} in
1951-52.

The only $1$-critical graph is $K_1$, and the only $2$-critical graph is $K_2$.
The only $3$-critical graphs are the odd cycles. For every $k\geq 4$ and every
$n\geq k+2$, there exists a $k$-critical $n$-vertex graph.
Let $f_k(n)$ be the minimum number of edges in a $k$-critical graph with $n$ vertices.
Since $\delta(G)\geq k-1$ for every $k$-critical $n$-vertex graph $G$,
\begin{equation}\label{in1}
f_k(n)\geq \frac{k-1}{2}n
\end{equation}
 for all $n\geq k$, $n\neq k+1$. Equality is achieved
for $n=k$ and for $k=3$ and $n$ odd. Brooks' Theorem \cite{Br} implies that for $k\geq 4$ and
$n\geq k+2$, the inequality in (\ref{in1}) is strict. In 1957, Dirac~\cite{D1} asked to determine
$f_k(n)$ and proved  that for $k\geq 4$ and
$n\geq k+2$,
\begin{equation}\label{in2}
f_k(n)\geq \frac{k-1}{2}n+\frac{k-3}{2}.
\end{equation}
The result is tight for $n=2k-1$ and yields $f_k(2k-1)=k^2-k-1$.
Dirac used his bound to evaluate chromatic number of graphs embedded into fixed surfaces.
Later, Kostochka and Stiebitz \cite{K2}  improved (\ref{in2}) to
\begin{equation}\label{Dirac Bound}
 f_k(n) \geq \frac{k-1}{2}n+k-3
\end{equation}
when $n \neq 2k - 1,k$.  This yields  $f_k(2k)=k^2-3$ and  $f_k(3k-2)=\frac{3k(k-1)}{2}-2$.
In his fundamental papers {\cite{G1,G2}},
Gallai found exact values of $f_k(n)$ for $k+2\leq n\leq 2k-1$:

\begin{theorem}[Gallai~\cite{G2}] \label{gallai1}
If $k \geq 4$ and $k+2\leq n\leq 2k-1$, then
$$f_k(n)=\frac{1}{2}  \left((k-1)n+(n-k)(2k-n)\right)-1.$$
 \end{theorem}
He also proved the following general bound for $k\geq 4$ and
$n\geq k+2$:
\begin{equation}\label{in3}
f_k(n)\geq \frac{k-1}{2}n+\frac{k-3}{2(k^2-3)}n.
\end{equation}
For large $n$, the bound is much stronger than bounds (\ref{in2}) and (\ref{Dirac Bound}).
Gallai in 1963 and Ore~\cite{O} in 1967 reiterated the question on finding $f_k(n)$.
Ore observed that Haj\'{o}s' construction implies
\begin{equation}\label{upper f_k}
f_k(n + k - 1) \leq f_k(n) + \frac{(k-2)(k+1)}2 = f_k(n) + (k-1)(k - \frac{2}{k-1})/2,
\end{equation}
which yields that $\phi_k:=\lim_{n\to \infty}\frac{f_k(n)}{n}$ exists and satisfies
\begin{equation}\label{in4}
\phi_k\leq \frac{k}{2}-\frac{1}{k-1}.
\end{equation}
Note that Gallai's bound gives $\phi_k\geq \frac{1}{2}\left(k-1+\frac{k-3}{k^2-3}\right)$.
Ore believed that Haj\'{o}s' construction was best possible.
\begin{conj} [Ore~\cite{O}] \label{Ore Conj} 
If $k \geq 4$, then 
$$ f_k(n + k - 1) = f_k(n) + (k-1)(k - \frac{2}{k-1})/2. $$
\end{conj}

Much later, Krivelevich~\cite{Kr2} improved Gallai's bound to
\begin{equation}\label{in5}
f_k(n)\geq \frac{k-1}{2}n+\frac{k-3}{2(k^2-2k-1)}n
\end{equation}
and demonstrated nice applications of his bound: he constructed graphs with high chromatic number and
low independence number such that the chromatic numbers of all their small subgraphs
are at most $3$ or $4$. We discuss a couple of his applications  in Subsection 6.3.
Then Kostochka and Stiebitz \cite{K2} proved that for $k\geq 6$ and $n\geq k+2$,
\begin{equation}\label{in6}
f_k(n)\geq \frac{k-1}{2}n+\frac{k-3}{k^2+6k-11-6/(k-2)}n.
\end{equation}

The problem of finding  $f_k(n)$ has attracted attention for more than $50$ years. It is Problem 5.3
in the monograph~\cite{J} and Problem~12 in the list of {\em $25$ pretty graph colouring problems}
by Jensen and Toft~\cite{J2}. It is one half of Problem P1 in~\cite[P. 347]{Tuza1}. Recently,
Farzad and Molloy~\cite{FM} have found the minimum number of edges in $4$-critical $n$-vertex
graphs in which the set of vertices of degree $3$ induces a connected subgraph.

The main result of the present paper is the following.

\begin{theorem} \label{k-critical}
If $k \geq 4$ and $G$ is $k$-critical, then
$ |E(G)|\geq  \left\lceil \frac{(k+1)(k-2)|V(G)|-k(k-3)}{2(k-1)}\right\rceil$. In other words, if
$k\geq 4$ and $n\geq k,\,n\neq k+1$, then
\begin{equation}\label{j20}
f_k(n)\geq  F(k,n):=\left\lceil \frac{(k+1)(k-2)n-k(k-3)}{2(k-1)}\right\rceil.
\end{equation}
\end{theorem}

This bound is exact for $k=4$ and every $n\geq 6$. For every $k\geq 5$, the bound is exact for every
$n \equiv 1\, ({\rm mod}\ k-1)$, $n\neq 1$. In particular, $\phi_k=\frac{k}{2}-\frac{1}{k-1}$ for every $k\geq 4$.
The result also confirms the above conjecture by Ore from 1967 for $k=4$ and every $n\geq 6$
and also for $k\geq 5$ and all
$n \equiv 1 \, ({\rm mod} \ k-1)$, $n\neq 1$. In the second half of the paper we derive some corollaries of the main result,
in particular, we give a very short proof of Gr\" otzsch' Theorem that every triangle-free planar graph is $3$-colorable.
Some further consequences are discussed in~\cite{BKY}.

Our proof of Theorem \ref{k-critical} is constructive.
This allows us to give an algorithm for coloring graphs with no dense subgraphs.
The idea of sparseness is expressed in terms of potentials.

\begin{defn}
For $R \subseteq V(G)$, define \emph{the $k$-potential of $R$} to be
\begin{equation}\label{rho}
\rho_{k,G}(R) = (k-2)(k+1)|R| - 2(k-1)|E(G[R])|.
\end{equation}
When there is no chance for confusion, we will use $\rho_k(R)$.
Let $P_k(G) = \min_{\emptyset \neq R \subseteq V(G)} \rho_k(R)$.
\end{defn}

\begin{theorem}\label{main balanced}
If $k \geq 4$, then every $n$-vertex graph $G$ with $P_k(G) > k(k-3)$ can be $(k-1)$-colored in $O(k^{3.5}n^{6.5}\log(n))$ time.
\end{theorem}

The restriction $P_k(G) > k(k-3)$ is sharp for every $k\geq 4$.
The next two corollaries follow from Theorems~\ref{k-critical} and~\ref{gallai1} and from   (\ref{upper f_k}).
Both will be proven in Section 5.

\begin{cor} \label{asymptotic}
For every $k \geq 4$ and $n\geq k+2$,
$$ 0 \leq f_k(n) - F(k,n) \leq \frac{k(k-1)}{8}-1.
$$
In particular,
$ \phi_k = \frac{k}2 - \frac1{k-1}.
$
\end{cor}

\begin{cor} \label{Ore cor}
For each fixed $k \geq 4$, Conjecture~\ref{Ore Conj} is true for all but at most
$\frac{k^3}{12} - \frac{k^2}{8}$
values of $n$.
\end{cor}

In Section 2 we prove several statements about list colorings that  will be used in our proofs.
In Section 3 we give definitions and prove several lemmas needed to prove Theorem~\ref{k-critical}
which will be proved 
in Section 4. %we prove Theorem \ref{k-critical}.
In Section 5 we discuss the sharpness of our result.
In Section 6 we present some applications.
In Section 7 we prove Theorem~\ref{main balanced}.
We finish the paper with some comments.

Our notation is standard. In particular, $\chi(G)$ denotes the chromatic number of graph $G$,
$G[W]$ is the subgraph of a graph or digraph $G$ induced by the vertex set $W$. For
a vertex $v$ in a graph $G$, $d_G(v)$ denotes the degree
of vertex $v$ in graph $G$, $N_G(v)$ is the set of neighbors of $v$ and $N_G[v]=N_G(v)\cup \{v\}$.
If the graph $G$ is clear from the context, we drop the subscript.

\section{Orientations and list colorings}
We consider loopless digraphs.
A {\em kernel} in a digraph $D$ is an independent set  $F$ of  vertices such that each vertex in $V(D)-F$
has an out-neighbor in $F$.

A digraph $D$ is {\em kernel-perfect} if for every $A\subseteq V(D)$, the digraph $D[A]$ has a
kernel. It is known that kernel-perfect orientations form a useful tool for list colorings.
Recall that a {\em list}  for a graph $G$ is a mapping $L$ of $V(G)$ into the family of finite
subsets of ${\bf N}$. For a given list $L$, a graph $G$ is $L$-{\em colorable}, if there exists
a coloring $f\,:\,V(G)\to {\bf N}$ such that $f(v)\in L(v)$ for every $v\in V(G)$ and
$f(v)\neq f(u)$ for every $uv\in E(G)$.
The following fact is well known but we include its proof for completeness.
\begin{lemma}[Folklore]\label{lem1} If $D$ is a kernel-perfect digraph and $L$ is a list such that
\begin{equation}\label{e1}
|L(v)|\geq 1+d^+(v)\qquad \mbox{ for every $v\in V(D)$,}
\end{equation}
 then $D$ is $L$-colorable.
\end{lemma}

{\bf Proof.} We use induction on $|V(D)|$. If $D$ has only one vertex, the statement is trivial.
Suppose the statement holds for all pairs $(D',L)$ satisfying (\ref{e1}) with $|V(D')|\leq n-1$.
Let $|V(D)|=n$ and $(D,L)$ satisfy (\ref{e1}). Let $v\in V(D)$ and $\alpha$ be a color present in
$L(v)$. Let $V_{\alpha}$ be the set of vertices  $x\in V(D)$  with $\alpha\in L(x)$.
Since $D$ is kernel-perfect, $D[V_{\alpha}]$ has a kernel $K$. Color all vertices of $K$ with
${\alpha}$ and consider $(D',L')$, where $D'=D-K$ and $L'(y)=L(y)-{\alpha}$ for all $y\in V(D')$.
Since the outdegree of every $x\in V_{\alpha}-K$ decreased by at least $1$, $(D',L')$ satisfies
(\ref{e1}), and so by the induction assumption has an $L'$-coloring. Together with coloring of $K$
 by ${\alpha}$, this yields an $L$-coloring of $D$, as claimed.\qed

It is known that every orientation of a bipartite multigraph is kernel-perfect. We prove a somewhat
stronger result.

\begin{lemma}\label{lem2} Let $A$ be an independent set in a graph $G$ and $B=V(G)-A$.
Let $D$ be the digraph obtained from $G$ by replacing each edge in $G[B]$ by a pair of
opposite arcs and by an arbitrary orientation of the edges connecting $A$ with $B$.
Then $D$ is kernel-perfect.
\end{lemma}

{\bf Proof.}  Let $D$ be a counter-example with the fewest vertices.
If every $b\in B$ has an outneighbor in $A$, then $A$ is a kernel. Otherwise,
 some $b\in B$ has no
outneighbors in $A$. Then $N(b)=N^-(b)$. We
 %include $b$ into a kernel and
 consider $D'=D-b-N^-(b)$. By the minimality of $D$, $D'$ has a kernel $K$. Then
$K+b$ is a kernel of $D$.\qed

For a graph $G$ and disjoint vertex subsets $A$ and $B$, let $G(A,B)$ denote the bipartite
graph with partite sets $A$ and $B$ whose edges are all edges of $G$ connecting $A$ with $B$.
The main result of this section is the following.

\begin{lemma}\label{lem3} Let $G$ be a $k$-critical graph. Let
 disjoint vertex subsets  $A$ and $B$ be such that\\
 (a) at least one of $A$ and $B$ is independent;\\
(b) $d(a)=k-1$ for every $a\in A$;\\
(c) $d(b)=k$ for every $b\in B$.\\
Then (i) $\delta(G(A,B))\leq 2$ and\\
(ii) either some $a\in A$ has at most one neighbor in $B$ or
some $b\in B$ has at most three neighbors in $A$.
\end{lemma}

{\bf Proof.} If $A\cup B=\emptyset$, then both statements are trivial. Otherwise,
since $G$ is
$k$-critical, there exists a $(k-1)$-coloring $f$ of $G-A-B$. Fix any such $f$.
For every $x\in A\cup B$, let $L(x)$ be the set of colors in $\{1,\ldots,k-1\}$
not used in $f$ on neighbors of $x$. Let $G'=G[A\cup B]$.
Then
\begin{equation}\label{e2}
\mbox{for every $a\in A$, $|L(a)|\geq d_{G'}(a)$, and
for every $b\in B$, $|L(b)|\geq d_{G'}(b)-1$.}
\end{equation}

CASE 1: $\delta(G(A,B))\geq 3$.
Let $G''$ be obtained from $G(A,B)$ by splitting each  $b\in B$ into
$\lceil d_{G(A,B)}(b)/3\rceil$ vertices of degree at most $3$. In particular,
a vertex $b$ of degree $3$ in $G(A,B)$ is not split. Graph $G''$ is bipartite
with partite sets $A$ and $B'$, where $B'$ is obtained from $B$.
The degree of each $a\in A$ in $G''$ is at least $3$, and the degree of each vertex
$b\in B'$ is at most $3$. So by Hall's Theorem, $G''$ has a matching $M$ covering $A$.
We construct digraph $D$ from $G'$ as follows:\\
(1)  replace each edge of $G[B]$ or in $G[A]$ (whichever is nonempty) with two opposite arcs,\\
(2) orient every edge of $G(A,B)$ corresponding to an edge in $M$ towards $A$,\\
(3) all other edges of $G(A,B)$ orient towards $B$.

By Lemma~\ref{lem2}, $D$ is kernel-perfect. Moreover, by (\ref{e2}), for every $a\in A$,
$d^+(a)=d_{G'}(a)-1\leq |L(a)|-1$, and for every $b\in B$,
$$d^+(b)\leq d_{G'}(b)-\lfloor \frac{2}{3}d_{G(A,B)}(b)\rfloor\leq (|L(b)|+1)-2=|L(b)|-1.
$$
Thus by Lemma~\ref{lem1}, $G'$ is $L$-colorable. But this means that $G$ is
$(k-1)$-colorable, a contradiction. This proves (i).

CASE 2: Each $a\in A$ has at least two neighbors in $B$ and each $b\in B$
has at least four neighbors in $A$. Then we obtain $G''$ by splitting each  $b\in B$ into
$\lceil d_{G(A,B)}(b)/2\rceil$ vertices of degree at most $2$.
Similarly to Case 1, graph $G''$ is bipartite
with partite sets $A$ and $B'$, where $B'$ is obtained from $B$.
The degree of each $a\in A$ in $G''$ is at least $2$, and the degree of each vertex
$b\in B'$ is at most $2$. So by Hall's Theorem, $G''$ has a matching $M$ covering $A$.
We construct digraph $D$ from $G''$ according to Rules (1)--(3) in Case 1.
Again, by Lemma~\ref{lem2}, $D$ is kernel-perfect, and by (\ref{e2}), for every $a\in A$,
$d^+(a)=d_{G'}(a)-1\leq |L(a)|-1$. For every $b\in B$, since $d_{G(A,B)}(b)\geq 4$, by (\ref{e2}),
 $$d^+(b)\leq d_{G''}(b)-\lfloor \frac{1}{2}d_{G(A,B)}(b)\rfloor\leq (|L(b)|+1)-2=|L(b)|-1.\qed$$

\begin{cor}\label{co1} Let $G$ be a $k$-critical graph. Let
 disjoint vertex subsets  $A$ and $B$ be such that\\
 (a) either $A$ or $B$ is independent;\\
(b) $d(a)=k-1$ for every $a\in A$;\\
(c) $d(b)=k$ for every $b\in B$;\\
(d) $|A|+|B|\geq 3$.\\
%Let $l=|A|$ and $h=|B|$.
 Then (i) $e(G(A,B))\leq 2(|A|+|B|)-4$ and
(ii) $e(G(A,B))\leq |A|+3|B|-3$.
\end{cor}

{\bf Proof.} First we prove (i) by induction on $|A|+|B|$. If $|A|+|B|=3$,
then since $G(A,B)$ is bipartite, it has at most $2=2\cdot 3-4$ edges.
Suppose now  that $|A|+|B|=m\geq 4$ and the corollary holds for $3\leq |A|+|B|\leq m-1$.
By Lemma~\ref{lem3}(i), $G(A,B)$ has a vertex $v$ of degree at most
two. By the minimality of $m$, $G(A,B)-v$ has at most $2(m-1)-4$ edges.
Then $e(G(A,B))\leq 2+2(m-1)-4=2m-4$, as claimed.

The base case $|A|+|B|=3$ for (ii) is slightly more complicated. If $|A|=3$, then
$e(G(A,B))=0=|A|+3|B|-3$. If $|B|\geq 1$, then $|A|+3|B|\geq 5$ and
$e(G(A,B))\leq 2=5-3\leq |A|+3|B|-3$.
The proof of the induction step is very similar to the previous paragraph, using Lemma~\ref{lem3}(ii).\qed

\section{Preliminary Results}
\begin{fact} \label{list of potentials} For the $k$-potential defined by (\ref{rho}), we have
\begin{enumerate}
	\item $\rho_{k,K_k}(V(K_k)) = k(k-3)$,
	\item $\rho_{k,K_1}(V(K_1)) = (k-2)(k+1)$,
	\item $\rho_{k,K_2}(V(K_2)) = 2(k^2-2k-1)$,	
	\item $\rho_{k,K_{k-1}}(V(K_{k-1})) = 2(k-2)(k-1)$.
\end{enumerate}
\end{fact}

A graph $H$ is {\em smaller than} graph $G$, if either $|E(G)|>|E(H)|$, or  $|E(G)|=|E(H)|$ and $G$ has fewer pairs of vertices with the same closed
neighborhood.
If $|V(G)| \geq |V(H)|$, $\rho_k(V(G)) \leq \rho_k(V(H))$, and equality does not hold in both cases, then $H$ is smaller than $G$.

Note that $(k - \frac2{k-1})|V(G)| > 2|E(G)| + \frac{k(k-3)}{k-1}$ is equivalent to $\rho_k(V(G)) > k(k-3)$.
Let $G$ be a minimal $k$-critical graph with respect to our relation "smaller" with $\rho_k(V(G)) > k(k-3)$.
This implies that 
\begin{equation}\label{j26}
\mbox{if $H$ is smaller than $G$ and $P_k(H) > k(k-3)$, then $H$ is $(k-1)$-colorable.}
\end{equation}

\begin{defn}\label{def1} For a graph $G$, a set $R\subset V(G)$ and a $(k-1)$-coloring $\phi$ of
$G[R]$, the graph $Y(G,R,\phi)$ is constructed as follows. First, for $i=1,\ldots,k-1$, let
$R'_i$ denote the set of vertices in $V(G)-R$ adjacent to at least one vertex  $v\in R$ with $\phi(v)=i$.
Second, let $X=\{x_1,\ldots,x_{k-1}\}$ be a set of new vertices disjoint from $V(G)$.
Now, let $Y=Y(G,R,\phi)$ be the graph with vertex set $ V( G) - R + X$, 
such that $Y[V(G)-R]=G-R$ and
$N(x_i) = R'_i \cup (\{x_1,\ldots,x_{k-1}\}-x_i)$ for $i=1,\ldots,k-1$.
\end{defn}

\begin{claim} \label{coloring Y}
Suppose $R \subset V(G)$ and $\phi$ is a $k-1$ coloring  of $G[R]$.
%If $Y(G, R, \phi)$ is smaller than $G$, t
Then $\chi(Y(G, R, \phi)) \geq k$.
\end{claim}
{\bf Proof.}
Let $G' = Y(G, R, \phi)$.
Suppose $G'$ has a $(k-1)$-coloring $\phi':V(G') \rightarrow C$.
By construction of $G'$, the colors of all $x_i$ in $\phi'$ are distinct.
By changing the names of the colors, we may assume that $\phi'(x_i) = i$ for $1 \leq i \leq k-1$.
By construction of $G'$, for all vertices $u \in R'_i$, $\phi'(u) \neq i$.
Therefore $\phi|_R \cup \phi'|_{V(G)-R}$ is a proper coloring of $G$,  a contradiction.
\qed

\begin{claim} \label{very small potential}
There is no $R \subsetneq V(G)$ with $|R| \geq 2$ and $\rho_{k,G}(R) \leq (k-2)(k+1)$.
\end{claim}
{\bf Proof.} Let $2\leq |R|<|V(G)|$ and
 $\rho_k(R)=m=\min\{\rho_k(W)\,:\, W\subsetneq V(G),\; |W|\geq 2\}$.
Suppose $m\leq (k-2)(k+1)$. Then $|R| \geq k$. Since $G$ is $k$-critical,
 $G[R]$ has a proper coloring $\phi:R \rightarrow C = \{1, \dots {k-1}\}$.
Let $G'=Y(G,R,\phi)$. By Claim~\ref{coloring Y},  $G'$ is not $(k-1)$-colorable.
Then it contains a $k$-critical subgraph $G''$. Let $W=V(G'')$. Since $|R| \geq k>|X|$
and $\rho_k(R)<\rho_k(X)$, $G''$ is smaller than $G$.
So, by the minimality
of $G$, $\rho_{k,G'}(W) \leq k(k-3)$.
Since $G$ is $k$-critical by itself, $W\cap X\neq \emptyset$.
Since every non-empty subset of $X$ has potential at least $(k-2)(k+1)$,
  $$\rho_{k,G}(W - X + R) \leq \rho_{k,G'}(W)-(k-2)(k+1)+m\leq m-2k+2.$$
Since $W-X+R\supset R$, $|W-X+R|\geq 2$. Since $\rho_{k,G}(W - X + R)<\rho_{k,G}( R)$,
by the choice of $R$, $W-X+R=V(G)$. But then $\rho_{k,G}(V(G)) \leq m-2k+2\leq k(k-3)$,
 a contradiction. \qed

\begin{lemma}\label{lem1 - k} Let $k-1\geq 2$ be an integer.
Let $R_*=\{u_1,\ldots,u_s\}$ be a vertex set and $w\,:\,R_* \to \{1,2,\ldots\}$ be an integral positive weight function on $R_*$
such that $w(u_1)+\ldots+w(u_s)\geq k-1$. Then for each $1\leq i\leq (k-1)/2$, there exists a
graph $H$ with $V(H)=R_*$ and
$|E(H)|=i$ such that for every independent set $M$ in $H$ with $|M|\geq 2$,
\begin{equation}\label{l01}
\mbox{ $\sum_{u\in R_*-M}w(u)\geq i$.}
\end{equation}
\end{lemma}

{\bf Proof.} We may assume that $w(u_1)\geq w(u_2)\geq\ldots\geq w(u_s)$.

CASE 1: $w(u_2)+\ldots+w(u_s)\leq i$. We let $E(H)=\{u_1u_j\,:\; 2\leq j\leq s\}$.
If $M$ is any independent set with $|M|\geq 2$,
then $u_1\notin M$ and witnesses that (\ref{lem1 - k}) holds.

CASE 2: $w(u_2)+\ldots+w(u_s)\geq i+1$. Choose the largest $j$ such that
$w(u_{j})+\ldots+w(u_s)\geq i$. Let $\alpha=i-w(u_{j+1})+\ldots+w(u_s)$.
 Since $ i\leq (k-1)/2$ and $w(u_1)+\ldots+w(u_s)\geq k-1$,
we also have $w(u_1)+\ldots+w(u_{j})\geq i + \alpha$. By the choice of $j$ and the ordering
of the vertices, $0 < \alpha\leq w(u_{j})\leq w(u_1)$.
We draw $\alpha$ edges connecting $u_1$ with $u_{j}$ and $i-\alpha$ edges
connecting $\{u_{j+1},\ldots,u_s\}$ with $\{u_1,\ldots,u_{j}\}$ so that
for each $\ell$, the degree of $u_\ell$ in the obtained multigraph $H$ is at most $w(u_\ell)$.
 Let $M$ be any nonempty independent set in $H$.
By the definition of $H$, since $M$ is independent,
$$\sum_{u\in R_*-M} w(u)\geq
\sum_{u\in R_*-M}d_H(u)\geq \frac{1}{2}\sum_{u\in R_*}d_H(u)=i,$$
as claimed. 
%Note that in this case, (\ref{l01}) holds for {\em every} independent set $M$,
%even if $|M|=1$.
If $H$ has multiple edges, we replace each set of multiple edges with a single edge.
\qed

\begin{claim} \label{small potential}
If $R \subsetneq V(G)$, $|R|\geq 2$ and $\rho_k(R) \leq 2(k-2)(k-1)$, then $R$ is a $K_{k-1}$.
\end{claim}
{\bf Proof.}
Let $R$ have the smallest $\rho_k(R)$ among $R \subsetneq V(G)$, $|R|\geq 2$.
Suppose $m=\rho_k(R) \leq 2(k-2)(k-1)$ and $G[R]\neq K_{k-1}$. Then $|R|\geq k$.
Let $i$ be the  integer such that
\begin{equation}\label{j14}
1+k(k-3)+2i(k-1)\leq \rho_k(R)\leq k(k-3)+2(i+1)(k-1).
\end{equation}
By Claim~\ref{very small potential}, $i\geq 1$.
Since for $k \geq 3$,
\begin{equation}\label{a22}
1+k(k-3)+\frac{k-1}{2}2(k-1)>2(k-2)(k-1),
\end{equation}
 we have $i\leq \frac{k-2}{2}$.

For $u\in R$, let $w(u)=|N(u)\cap (V(G)-R)|$.
Let $R_* = \{u \in R\,:\, w(u)\geq 1\}$.
Because $\kappa(G) \geq 2$, $|R_*| \geq 2$.
Since $G$ is $k$-critical, $\sum_{u\in R_*}w(u)=|E_G(R,V(G)-R)|\geq k-1$.
So by Lemma~\ref{lem1 - k}, we can add to $G[R_*]$ a set $E_0$ of at most $i$ edges so that for every independent subset $M$ of $R_*$ in $G\cup E_0$ with $|M|\geq 2$, (\ref{l01}) holds.
Let $H=G[R]\cup E_0$.
Note that $|E(G)| - |E(G[R])| \geq k-1 > i$, so $H$ is smaller than $G$.
By the minimality of $\rho_k(R)$ and the definition of $i$, for every $U\subseteq R$ with $|U|\geq 2$,
$$\rho_{k,H}(U)\geq \rho_{k,G}(U)-2i(k-1)\geq \rho_{k,G}(R)-2i(k-1) \geq 1+k(k-3).$$
Thus $P_k(H)\geq 1+k(k-3)$, and by (\ref{j26}) $H$ has
 a proper $(k-1)$-coloring $\phi$ with colors in $C= \{1, \ldots,{k-1}\}$.

As in Claim \ref{very small potential}, we let $G'=Y(G,R,\phi)$.
Since $|R|\geq k$, $|V(G')|<|V(G)|$. Since
$$\rho_{k,G'}(V(G'))=\rho_{k,G}(V(G))-\rho_k(R)+\rho_k(X)\geq \rho_{k,G}(V(G)),$$
$|E(G')|<|E(G)|$ and so $G'$ is smaller than $G$.
By Claim \ref{coloring Y}, $G'$ is not  $(k-1)$-colorable. Thus
$G'$ contains  a $k$-critical subgraph $G''$. Let $W=V(G'')$. By the minimality
of $G$, $\rho_{k,G'}(W) \leq k(k-3)$.
Since $G$ is $k$-critical by itself, $W\cap X\neq \emptyset$.

Since every subset of $X$ with at least two vertices has potential at least $2(k-2)(k-1)$,
 if $|W \cap X| \geq 2$ then  $\rho_{k,G}(W - X + R) \leq \rho_{k,G'}(W)\leq k(k-3)$,  a contradiction again.
So, without loss of generality, assume that $X \cap W = \{x_1\}$.
But then
\begin{equation}\label{e2 - k}
\rho_{k,G}(W - \{x_1\} + R)  \leq  (\rho_{k,G'}(W)-(k-2)(k+1)) + \rho_{k,G}(R)\leq    \rho_{k,G}(R) - 2k + 2.
\end{equation}
%$$ \leq    \rho_{k,G}(R) - 2k + 2.$$
By the minimality of $\rho_{k,G}(R)$, $W - \{x_1\} + R = V(G)$.
This implies that $W = V(G') - X+x_1$.

Let $R_1 = \{u \in R_*: \phi(u) = \phi(x_1)\}$.
If $|R_1|=1$, then
$$\rho_{k,G}(W-x_1 \cup R_1)= \rho_{k,H}(W) \leq k(k-3),$$
a contradiction.
Thus, $|R_1|\geq 2$.
Since $R_1$ is an independent set,
 by the construction of $H$, at least $i$ edges connect the vertices in $R_*-R_1$ with $V(G)-R$.
These edges were not accounted in (\ref{e2 - k}).
So, in this case  instead of (\ref{e2 - k}), we have
\begin{eqnarray*}
	\rho_{k,G}(W - \{x_1\} + R) & \leq & \rho_{k,G'}(W)-(k-2)(k+1)-2i(k-1) + \rho_{k,G}(R) \\
	 & \leq & \rho_{k,G}(R)-2k + 2-2i(k-1) \\
	 & = & \rho_{k,G}(R) - 2(i+1)(k-1) \\
	 & \leq & k(k-3), \\
\end{eqnarray*}
%which is
a contradiction.
\qed

\begin{claim} \label{clusters in cliques}
If $d(x) = d(y) = k-1$ and $x$ and $y$ are in the same $(k-1)$-clique, then $N[x] = N[y]$.
\end{claim}
{\bf Proof.}
By contradiction, assume that $d(x_1) = d(x_2) = k-1$, $N(x_1) =X-x_1+a$,  $N(x_2) =X-x_2+b$, and  $a \neq b$.
If $ab \in E(G)$, then define $G' = G-\{x_1,x_2\}$.
Otherwise define $G' = G - \{x_1,x_2\} + ab$.
Because $\rho_{k,G}(W) \geq 2(k-2)(k-1)$ for all $W \subseteq G-\{x_1,x_2\}$ with $|W| \geq 2$, and adding an edge decreases the potential of a set by $2(k-1)$,
$$P_k(G') \geq \min\{(k-2)(k+1),2(k-2)(k-1)-2(k-1)\}>1+k(k-3).$$
So, since $G'$ cannot contain $k$-critical subgraphs, it
%It follows that $G'$
 has a proper $(k-1)$-coloring $\phi'$ with $\phi'(a) \neq \phi'(b)$.
This easily extends to a proper $(k-1)$-coloring of $V(G)$.
\qed

\begin{defn}
A { \em cluster} is a maximal set $R \subseteq V(G)$ such that for every $x \in R$, $d(x) = k-1$ and for every pair $x,y \in R$, $N[x] = N[y]$.
\end{defn}

\begin{claim} \label{small clusters}
Let $C$ be a cluster.
Then $|C| \leq k-3$.
Furthermore, if $C$ is in a $(k-1)$-clique $X$, then $|C| \leq \frac{k-1}{2}$.
\end{claim}
{\bf Proof.}
A cluster with $k-2$ vertices plus its two neighbors would form a set of potential
at most $k(k-3) + 2(k-1)$, which is less than $2(k-2)(k-1)$ when $k \geq 4$.

Let $\{v\} = N(C) - X$.
If $|C| \geq \left\lceil k/2\right\rceil$, then $\rho_k(X + v) \leq 2(k-2)(k-1) - 2$,  a contradiction.
\qed

\begin{claim} \label{adjacent k-1}
Let $xy \in E(G)$, $N[x] \neq N[y]$, $x$ is in a cluster of size $s$, $y$ is in a cluster of size $t$, and $s \geq t$.
Then $x$ is in a $(k-1)$-clique.
Furthermore, $t = 1$.
\end{claim}
{\bf Proof.}
Assume that $x$ is not in a $(k-1)$-clique.
Let $G' = G - y + x'$, where $N[x'] = N[x]$.
We have $|E(G')| = |E(G)|$. If two vertices $z$ and $z'$ distinct from  $y$ had the same
closed neighborhood in $G$, then they  also have the same closed neighborhood in $G'$. Thus, since
 the cluster containing $x$ is at least
as large as the one containing $y$, $G'$ is smaller than $G$ in our ordering.
If  $G'$ has a $(k-1)$-coloring $\phi':V(G') \rightarrow C = \{1, 2, \dots {k-1}\}$, then
we extend it to a proper $(k-1)$-coloring $\phi$ of $G$ as follows:
define $\phi|_{V(G)-x-y} = \phi'|_{V(G')-x-x'}$, then  choose $\phi(y) \in C - (\phi'(N(y) - x))$,
and $\phi(x) \in \{\phi'(x), \phi'(x')\} - \{\phi(y)\}$.

So, $\chi(G')\geq k$ and $G'$ contains a $k$-critical subgraph $G''$. Let $W=V(G'')$.
Since $G''$ is smaller than $G$, $\rho_{k, G'}(W) \leq k(k-3)$. Since $G''$ is not a subgraph of $G$,
 $x' \in W$.
Then $\rho_{k,G}(W - x') \leq k(k-3) - (k-2)(k+1) + 2(k-1)(k-1) = 2(k-2)(k-1)$.
This contradicts Claim \ref{small potential} because $y \notin W-x'$ and so $W-x'\neq V(G)$.
\qed

\section{Proof of Theorem \ref{k-critical}}

\subsection{Case $k = 4$}
%We define $\rho(R) = \rho_4(R) / 2 = 5|R| - 3|E(G[R])|$ and $P(G) = P_4(G)/2$.
%Let $G$ be still the minimum counterexample with respect to the order defined above.
%The theorem is then equivalent to

%\begin{theorem}\label{main 3 balanced}
% $\rho_k(V(G)) \leq 2$.
%\end{theorem}

\begin{claim} \label{1 K3}
Each edge of $G$ is in at most $1$ triangle.
Moreover, each cluster has only one vertex.
\end{claim}
{\bf Proof.} The vertex set of a subgraph with
$4$ vertices and $5$ edges has  potential $10$, which contradicts Claim~\ref{small potential}.
A cluster of size two would create an edge shared by two triangles.
\qed

\begin{claim} \label{two 3 etc}
Each vertex with degree $3$ has at most $1$ neighbor with degree $3$.
\end{claim}
{\bf Proof.}
This follows directly from Claims \ref{1 K3} and \ref{adjacent k-1}.
\qed

We will now use discharging to show that $|E(G)| \geq \frac{5}3 |V(G)|$, which will finish the proof to
the case $k=4$.
Each vertex begins with charge equal to its degree.
If $d(v) \geq 4$, then $v$ gives charge $\frac16$ to each neighbor with degree $3$.
Note that $v$ will be left with charge at least $\frac56 d(v) \geq \frac{10}{3}$.
By Claim \ref{two 3 etc}, each vertex of degree $3$ will end with charge at least $3 + \frac26 = \frac{10}3$.
\qed

\subsection{Case $k = 5$}
%We define $\rho'(R) = \rho_5(R) / 2 = 9 |R| - 4 |E(G[R])|$ and $P'(G) = P_5(G)/2$.
%Let $G$ be the minimum counterexample with respect to the order defined above.
%The theorem is then equivalent to

%\begin{theorem}\label{main 4 balanced}
 %$P'(G) \leq 5$.
%\end{theorem}

\begin{claim} \label{4 no neighbors}
Each cluster has only one vertex.
\end{claim}
{\bf Proof.}
Assume $N[x] = N[y]$ and $d(x) = d(y) = 4$.
Because $G$ does not contain a $K_5$, there exist $a,b \in N[x]$ such that $ab \notin E(G)$.
We obtain $G'$ from $G$ by deleting $x$ and $y$  and gluing $a$ with $b$.
If $G'$ is $4$-colorable, then so is $G$.
This is because a $4$-coloring of $G'$ will have at most $2$ colors on $N[x] - \{x,y\}$ 
and therefore could be extended greedily to $x$ and $y$.

So $G'$ contains a $k$-critical subgraph $G''$. Let $W=V(G'')$. Since $G''$ is smaller than $G$,
 $\rho_{5,G'}(W) \leq 10$. Since $G''$ is not a subgraph of $G$, $a*b \in W$.
But then $\rho_{5,G}(W - a*b + a + b + x + y) \leq 10 + 54 - 40 = 24$.
Because $ab \notin E(G)$, $W - a*b + a + b + x + y$ is not a $K_4$.
By Claim \ref{small potential}, $W - a*b + a + b + x + y = V(G)$.
But then we did not account for two of the edges incident to $\{x,y\}$, so $\rho'_{G}(W - a*b + a + b + x + y) \leq
24-2\cdot 8=8$,  a contradiction.
\qed

\begin{claim}\label{c23}
 Each $K_4$-subgraph of $G$ contains at most one vertex  with degree $4$.
If $d(x) = d(y) = 4$ and $xy \in E(G)$, then each of $x$ and $y$ is in a $K_4$.
\end{claim}
{\bf Proof.}
The first statement follows from Claims \ref{clusters in cliques} and \ref{4 no neighbors}.
The second statement follows from Claims \ref{adjacent k-1} and \ref{4 no neighbors}.
\qed

\begin{defn}
We define $H \subseteq V(G)$ to be the set of vertices of degree $5$
not in a $K_4$, and $L \subseteq V(G)$ to be the set of vertices of degree $4$ not in a $K_4$.
Set $\ell=|L|$, $h=|H|$ and $e_0 = |E(L, H)|$.
\end{defn}

\begin{claim} \label{HL charge}
$e_0 \leq  3h + \ell$.
\end{claim}
{\bf Proof.}
This is trivial if $h+\ell\leq 2$ and follows from
 Corollary~\ref{co1}(ii) and Claim~\ref{c23} for $h+\ell\geq 3$.
\qed

We will do discharging in two stages.
Let every vertex  $v\in V(G)$ have initial charge $d(v)$.
The first half of discharging has one rule:

{\bf Rule R1:} Each vertex in $V(G)-H$ with degree at least $5$ gives charge $1/6$ to each neighbor.

%{\bf Rule R2:} Each vertex in $H$ gives $1/9$ to each neighbor in $L$.

\begin{claim} \label{5 discharging 1}
After the first round of discharging, each vertex in $V(G) - H - L$ has charge at least $4.5$.
\end{claim}
{\bf Proof.}
Let $v \in V(G) - H - L$.
If $d(v) = 4$, then $v$ receives $1/6$ from at least $3$ neighbors and gives no charge.
If $d(v) = 5$, then $v$ gives $1/6$ to $5$ neighbors, but receives $1/6$ from at least $2$ neighbors.
If $d(v) \geq 6$, then $v$ is left with charge at least $5d(v)/6 \geq 4.5$.
\qed

For the second round of discharging, all charge in $H \cup L$
is taken up and distributed evenly among the vertices in $H \cup L$.

\begin{claim} \label{5 discharging 2}
After the first round of discharging, the sum of the charges on the vertices in $H \cup L$ is at least $4.5|H \cup L|$.
\end{claim}
{\bf Proof.} By Rule R1, vertices in $L$ receive from outside of $H\cup L$ the charge at least $\frac{1}{6}(4\ell-|E(H,L)|)$.
By Claim \ref{HL charge}, $|E(H,L)| \leq  3h + \ell $. So, the total charge on $H\cup L$ is
at least
$$5h+4\ell+\frac{1}{6}(4\ell-(3h+\ell))=4.5(h+\ell),$$
as claimed.
\qed

Combining Claims \ref{5 discharging 1} and \ref{5 discharging 2},
the average degree of  the vertices in $G$ is at least $4.5$, a contradiction.

\subsection{Case $k \geq 6$}

\begin{claim} \label{big neighbors}
Let $T$ be a cluster in $G$ and $t=|T| \geq 2$.\\
(a) If $N(T)\cup T$ does not contain $K_{k-1}$, then $d_G(v)\geq k-1+t$ for every $v\in N(T)-T$;\\
(b) If $N(T)\cup T$  contains a $K_{k-1}$ with vertex set $X$, then $d_G(v)\geq k-1+t$ for every $v\in X-T$.\end{claim}
{\bf Proof.} Let $v\in N(T)-T$ such that $k \leq d(v) \leq k-2 + t$ and if
$N(T)\cup T$  contains a $K_{k-1}$ with vertex set $X$, then $v\in X$.
Since $\rho_{k,G}(N(T) \cup T) > (k+1)(k-2)$, $T$ is contained in at most one $(k-1)$-clique,
and so
\begin{equation}\label{j15}
\mbox{ $N(T)\cup T-v$ does not contain $K_{k-1}$.}
\end{equation}
%Assume that $k \leq d(v) \leq k-2 + t$ and $u \in T$.
By the choice of $v$, $|N(v) - T| \leq k-2$.
Let  $u \in T$ and $G' = G - v + u'$, where $N[u'] = N[u]$.
Suppose $G'$ has a $(k-1)$-coloring $\phi':V(G') \rightarrow C =\{1, \dots {k-1}\}$.
Then there is a $(k-1)$-coloring $\phi$ of $G$ as follows: set $\phi|_{V(G) - T - v} = \phi'|_{V(G') - T - u'}$, $\phi(v) \in C - \phi'(N(v) - T)$, and then color $T$ using colors from $\phi'(T \cup u') - \phi(v)$.
This is a contradiction, so there is no $(k-1)$-coloring of $G'$.
Thus $G'$ contains a $k$-critical subgraph $G''$. Let $W=V(G'')$.

Because $d_G(v) \geq k$ and $d_{G'}(u') = k-1$, $|E(G')| < |E(G)|$.
So, $G''$ is smaller than $G$ and hence
 $\rho_{k,G'}(W) \leq k(k-3)$. Since $G''$ is not a subgraph of $G$,
 $u' \in W$.
By symmetry, it follows that $T \subset W$.
But then
$$ \rho_{k,G}(W-u') \leq k(k-3) - (k-2)(k+1) + 2(k-1)(k-1) = 2(k-2)(k-1).$$
This implies that $G[W - u']$ is a $K_{k-1}$, a contradiction to (\ref{j15}).
\qed

\begin{claim} \label{cliques have k+1}
Suppose $v$ is the unique vertex with degree $k-1$ in a $(k-1)$-clique $X$.
Then $X$ contains at least $(k-1)/2$ vertices with degree at least $k+1$.
\end{claim}
{\bf Proof.}
Let $\{u\} = N(v) - X$.
Assume that $X$ contains at least $k/2-1$ vertices with degree $k$.
Note that $|N(u) \cap X| < k/2$, so there exists a $w \in X$ such that $uw \notin E(G)$ and $d(w) \leq k$.
Let $N(w) - X = \{a,b\}$.
Let $G'$  be obtained from $G-v$ by adding edges $ua$ and $ub$.  %if they do not already exist.

If $G'$ is not $(k-1)$-colorable, then it  contains a $k$-critical subgraph $G''$. Let $W=V(G'')$.
Since $|E(G')| < |E(G)|$, $G''$ is smaller than $G$ and so, $\rho_{k,G'}(W) \leq k(k-3)$.
If $W = V(G')$, then $\rho_{k,G}(V(G)) \leq k(k-3) + (k-2)(k+1)(1) - 2(k-1)(k-3) < k(k-3)$ when $k \geq 6$.
If $W \neq V(G')$ then $\rho_{k,G}(W) \leq k(k-3) + 2(k-1)(2) < 2(k-2)(k-1)$,  a contradiction.

Thus $G'$ has a $(k-1)$-coloring $f$. If $f(u)$ is not used
on $X-w-v$, then we  recolor $w$ with $f(u)$. So, anyway $v$ will have two neighbors of color $f(u)$,
and we can extend the  $(k-1)$-coloring to $v$.
\qed

\begin{claim} \label{k=6}
If $k=6$ and  a cluster $C$  is contained in a $5$-clique $X$, then $|C| = 1$.
\end{claim}
{\bf Proof.}
By Claim \ref{small clusters}, assume that $C = \{v_1, v_2\}$.
Let $N(v_1) - X = \{y\}$ and $\{u, u', u''\} = X-C$.
Obtain $G'$ from $G-C$ by gluing $u$  to $y$.

Suppose that $G'$ has a $5$-coloring.
We will extend this coloring to a coloring on $G$ by greedily assigning colors to $C$.
This can be done because only $3$ different colors appear on the vertices $\{u, u', u'', y\}$.
So we may assume that $\chi(G') \geq 6$. Then $G'$ contains a $k$-critical subgraph $G''$. Let $W=V(G'')$.
Because $|E(G')| < |E(G)|$, $\rho_{6,G'}(W) \leq 18$. Since $G''$ is not a subgraph of $G$,
 $u*y \in W$.
Let $t = | \{u', u''\} \cap W|$.

{\it Case 1:} $t=0$.
Then $\rho_{6, G}(W-u*y + y + X) \leq 18 + 28(5) - 10(12) = 38$.
By Claim \ref{small potential}, $W-u*y + y + X = V(G)$.
But then we did not account for edges in $E(\{u', u''\}, V(G)-X)$.
Thus $\rho_{6, G}(V(G)) \leq 38-2\cdot 10=18$.

{\it Case 2:} $t = 1$.
Then $\rho_{6, G}(W-u*y + y + u + C) \leq 18 + 28(3) - 10(7) = 32$.
This is a contradiction to Claim \ref{small potential}
because $V(G) \neq \left( W-u*y + y + u + C \right)$.

{\it Case 3:} $t = 2$.
Then $\rho_{6, G}(W-u*y + y + u + C) \leq 18 + 28(3) - 10(9) = 12$, which is a contradiction.
\qed

\begin{defn} \label{L,H,e0}
We partition $V(G)$ into four classes: $L_0$, $L_1$, $H_0$, and $H_1$.
Let $H_0$ be the set of vertices with degree $k$, $H_1$ be the set of vertices with degree at least $k+1$, and $H = H_0 \cup H_1$.
Let
$$L = \{u \in V(G): d(u) = k-1\},
$$
$$L_0 = \{u \in L : N(u) \subseteq H\},
$$
and
$$ L_1 = L - L_0.
$$
Set $\ell=|L_0|$, $h=|H_0|$ and $e_0 = |E(L_0, H_0)|$.
\end{defn}

\begin{claim} \label{G_0 charge}
$e_0 \leq 2(\ell+h)$.
\end{claim}
{\bf Proof.} 
This is trivial if $h+\ell\leq 2$ and follows from
 Corollary~\ref{co1}(i)  for $h+\ell\geq 3$.
\qed

Let every vertex  $v\in V(G)$ have initial charge $d(v)$.
We first do a half-discharging with two rules:

{\bf Rule R1:} Each vertex in $H_1$ keeps for itself charge $k-2/(k-1)$ and distributes the rest equally
among its neighbors of degree $k-1$.

{\bf Rule R2:} If a $K_{k-1}$-subgraph $C$ contains $s$ $(k-1)$-vertices adjacent to a $(k-1)$-vertex $x$
outside of $C$ and not in a $K_{k-1}$, then each of these $s$  vertices gives charge $\frac{k-3}{s(k-1)}$ to $x$.

\begin{claim}
Each vertex in $H_1$ donates at least $\frac{1}{k-1}$ charge to each neighbor of degree $k-1$.
\end{claim}
{\bf Proof.}
If $v \in H_1$, then $v$ donates at least $\frac{d(v) - k + 2/(k-1)}{d(v)}$ to each neighbor.
Note that this function increases as $d(v)$ increases, so the charge is minimized when $d(v) = k+1$.
But then each vertex gets charge at least $(1 + 2/(k-1))/(k+1)=1/(k-1)$.
\qed

\begin{claim}
Each vertex in $L_1$ has charge at least $k-2/(k-1)$.
\end{claim}
{\bf Proof.}
Let $v \in L_1$ be in a cluster $C$ of size $t$.

{\it Case 1:} $v$ is in a $(k-1)$-clique $X$ and $t \geq 2$.
By Claim \ref{k=6}, this case only applies when $k \geq 7$.

By Claim \ref{big neighbors} each vertex in $X-C$ has degree at least $k - 1 + t \geq k+1$, and therefore $X-C \subseteq H_1$.
Furthermore, each vertex in $X-C$ has at least $k-2-t$ neighbors with degree at least $k$.
Therefore each vertex $u \in (X-C)$ donates charge at least $\frac{d(u) - k + 2/(k-1)}{d(u)-k+2+t}$  to each neighbor of degree $k-1$.
Note that this function increases as $d(u)$ increases, so the charge is minimized when $d(u) = k-1+t$.
It follows that $u$ gives to $v$  charge at least $\frac{t-1+2/(k-1)}{2t+1}$.

So, $v$ has charge at least $k-1 + (k-1-t)(\frac{t-1 + 2/(k-1)}{2t+1}) - \frac{k-3}{t(k-1)}$, which we claim is at least $k-2/(k-1)$.
Let
$$ g_1(t) =(k-1-t)((t-1)(k-1)+2) - (2t+1)(k-3)(1+\frac{1}{t}).
$$
We claim that $g_1(t) \geq 0$, which is equivalent to $v$ having charge at least $k - 2/(k-1)$.
Let
$$ \widetilde{g}_1(t) =(k-1-t)((t-1)(k-1)+2) - (2t+1)(k-3)(3/2).
$$
Note that $\widetilde{g}_1(t) \leq g_1(t)$ when $t \geq 2$, so we need to show that $\widetilde{g}_1(t) \geq 0$ on the appropriate domain.
$\widetilde{g}_1(t)$ is quadratic with a negative coefficient at $t^2$, so it suffices to check its values at the boundaries.
They are
$$ \widetilde{g}_1(2)  = (k-3)(k-6.5)
$$
and
\begin{eqnarray*}
  4 \widetilde{g}_1(\frac{k-1}{2}) & = & (k-1)\left( (k-3)(k-1) + 4 \right) - 6 k (k-3) \\
  		& = & k^3 - 11 k^2 + 29 k - 7 \\
  		& = &  (k-7)(k^2 - 4k + 1).
\end{eqnarray*}
Each of these values is non-negative when $k \geq 7$.

{\it Case 2:} $t \geq 2$ and $v$ is not in a $(k-1)$-clique.
By Claim \ref{big neighbors}, each neighbor of $v$ outside of $C$ has degree at least $k - 1 + t \geq k+1$ and is in $H_1$.
Therefore $v$ has charge at least $k-1 + (k-t)(\frac{t-1 + 2/(k-1)}{k-1+t})$.
We define
\begin{eqnarray*}
	g_2(t) &=& (k-t)(t-1+\frac{2}{k-1}) - \frac{k-3}{k-1}(k-1+t)\\
		&=& t(k-t) - 2(1 - \frac{2}{k-1})(k-1)\\
		& = & t(k-t) - 2(k - 3).
\end{eqnarray*}
Note that $g_2(t) \geq 0$ is equivalent to $v$ having charge at least $k - 2/(k-1)$.
The function $g_2(t)$ is quadratic with a negative coefficient at $t^2$, so it suffices to check
its values at the boundaries.
They are
$$ g_2(2) = 2(k-2) - 2(k-3) = 2
$$
and
$$ g_2(k-3) = (k-3)(3) - 2(k-3) = k-3.
$$
Each of these values is positive.

{\it Case 3:} $t = 1$.
If $v$ is not in a $(k-1)$-clique $X$, then by Claim \ref{adjacent k-1} the vertex adjacent to $v$ with degree $k-1$ is in a $(k-1)$-clique and cluster of size at least $2$.
In this case $v$ will recieve charge $(k-3)/(k-1)$ in total from that cluster.
Therefore we may assume that $v$ is in a $(k-1)$-clique $X$.

By Claim \ref{cliques have k+1}, there exists a $Y \subset X$ such that $|Y| \geq \frac{k-1}2$ and every vertex in $Y$ has degree at least $k+1$.
Furthermore, each vertex in $Y$ has at least $k-3$ neighbors with degree at least $k$.
Therefore each vertex $u \in Y$ donates at least $\frac{d(u) - k + 2/(k-1)}{d(u)-k+3}$ charge to each neighbor of degree $k-1$.
Note that this function increases as $d(u)$ increases, so the charge is minimized when $d(u) = k+1$.
It follows that $u$ gives to $v$ charge at least $\frac{1+2/(k-1)}{4}$, and $v$ has charge at least
$$k-1 + \frac{k-1}2\left(\frac{1 + 2/(k-1)}{4}\right) = k + \frac{k-7}{8},
$$
which is at least $k-2/(k-1)$ when $k \geq 6$.
\qed

We then observe that after that half-discharging,\\
a) the charge of each vertex in $H_1\cup L_1$ is at least $k-2/(k-1)$;\\
b) the charges of vertices in $H_0$ did not decrease;\\
c) along every edge from $H_1$ to $L_0$ the charge at least $1/(k-1)$ is sent.

Thus by Claim \ref{G_0 charge}, the total charge $F$ of the vertices in $H_0\cup L_0$ is at least
$$kh+(k-1)\ell+\frac{1}{k-1}\left(\ell(k-1)-e(G_0)\right)\geq k(h+\ell)-\frac{1}{k-1}2(h+\ell)=
(h+\ell)\left(k-\frac{2}{k-1}\right),
$$
and so by a), the total charge of all the vertices of $G$
is at least $n\left(k-\frac{2}{k-1}\right)$, a contradiction.\qed

\section{Sharpness}
%{Minimum number of edges in a $k$-critical graph}
%\section{Theorem \ref{k-critical}}

%\section{Minimum number of edges in a $k$-critical graph}
%By Theorem \ref{k-critical},
%When $k \geq 4$, we have
%\begin{equation}\label{lower f_k}
%f_k(n) \geq \left\lceil\frac12\left( (k - \frac2{k-1})n - \frac{k(k-3)}{k-1} \right)\right\rceil.
%\end{equation}

The next statement shows some cases when the bound (\ref{j20}) of Theorem~\ref{k-critical}
is exact.
\begin{theorem}\label{th4}
If one of the following holds:
\begin{enumerate}
	\item $n \equiv 1\, ({\rm mod}\ k-1)$ and $n \geq k$,
	\item $k = 4$, $n \neq 5$, and $n \geq 4$, or
	\item $k = 5$, $n \equiv 2\, ({\rm mod}\ 4)$, and $n \geq 10$,
\end{enumerate}
then
$$ f_k(n) =  F(k,n) = \left\lceil \frac12\left( (k - \frac2{k-1})n - \frac{k(k-3)}{k-1} \right) \right\rceil.
$$
\end{theorem}
{\bf Proof.}
By (\ref{upper f_k}), we only need to show that (\ref{j20}) is tight when
\begin{enumerate}
	\item $n = k$,
	\item $k = 4$, $n = 6$,
	\item $k = 4$, $n = 8$, and
	\item $k = 5$, $n = 10$.
\end{enumerate}

\begin{figure}[htbp]
\begin{center}
\includegraphics[height=4cm]{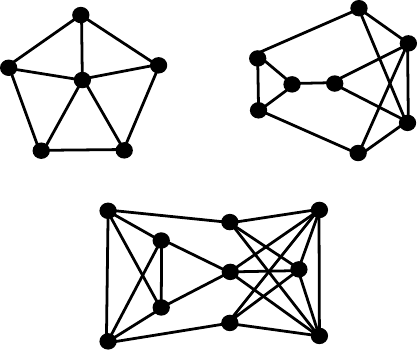}

\caption{Minimal $k$-critical graphs.}
\label{fk examples}
\end{center}
\end{figure}

The first case follows from $K_k$.
The other three cases follow from Figure \ref{fk examples}.
\qed

By Theorem~\ref{gallai1}, (\ref{j20}) is not sharp when $k\geq 5$ and $k+2\leq n\leq 2k-2$.
Probably, (\ref{j20}) is not sharp in case  of  $n$ not covered by Theorem~\ref{th4}.

Now we prove Corollary~\ref{asymptotic}. First, we restate it:

\bigskip\noindent
{\bf Corollary~\ref{asymptotic}}
{\em For $k \geq 4$,
$ 0 \leq f_k(n) - F(k,n) \leq (1+o(1))\frac{k^2}{8}.
$
In particular,
$ \phi_k = \frac{k}2 - \frac1{k-1}.
$}\\

\noindent
{\bf Proof.} By Theorem~\ref{th4}, the corollary holds for $k=4$. Let $k\geq 5$.
By  (\ref{upper f_k}) and Theorem~\ref{k-critical}, for every $n\geq k,\, n\neq k+1$,
$$f_k(n+(k-1))-F(k,n+(k-1))\leq f_k(n)-F(k,n).$$ Thus, it is enough to check the inequality for
$k+2\leq n\leq 2k$. There exists a $k$-critical $2k$-vertex graph with $k^2-3$ edges.
So,
$$f_k(2k)-F(k,2k)\leq k^2-3-\frac{(k+1)(k-2)2k-k(k-3)}{2(k-1)}\leq \frac{k(k-3)}{2(k-1)}<\frac{k-2}{2},$$
and by the integrality of $f_k$ and $F$,  $f_k(2k)-F(k,2k)\leq \frac{k-3}{2}$.

By Theorems~\ref{k-critical} and~\ref{gallai1}, for $k+2\leq n\leq 2k-1$,
\begin{equation} \label{separation}
f_k(n)-F(k,n)\leq \left(\frac{1}{2}  \left((k-1)n+(n-k)(2k-n)\right)-1\right)-\frac{(k+1)(k-2)n-k(k-3)}{2(k-1)}
\end{equation}
$$=-1+\frac{1}{2}\left[(n-k)\left(2k-\frac{k-3}{k-1}-n\right)\right].$$
For every fixed $k$, the maximum of the last expression (quadratic in $n$) is attained at
$n=\frac{1}{2}\left(k+2k-\frac{k-3}{k-1}\right)$. If $k\geq 5$, then the closest half-integer to this point
is $\frac{3k-1}{2}$. Thus,
$$f_k(n)-F(k,n)\leq f_k(\frac{3k-1}{2})-F(k,\frac{3k-1}{2})\leq -1+\frac{1}{2}\left[\frac{k-1}{2}\left(\frac{k+1}{2}
-\frac{k-3}{k-1}\right)\right]
$$
$$< -1+\frac{k-1}{4}\frac{k}{2}=-1+\frac{k(k-1)}{8}.\qed$$

\medskip
In particular, by the integrality of $f_k$ and $F$, $f_5(n)-F(5,n)\leq 1$ for all $n\geq 7$.

Now we prove Corollary~\ref{Ore cor}. First, we restate it:

\bigskip\noindent
{\bf Corollary~\ref{Ore cor}}
{\em If $k \geq 4$, then for all but $ \frac{k^3}{12} - \frac{k^2}{8}$ values of $n\geq k+2$, 
$$ f_k(n + k - 1) = f_k(n) + (k-1)(k - \frac{2}{k-1})/2. $$
}\\

\noindent
{\bf Proof.} By Theorem~\ref{th4}, the corollary holds for $k=4$. 
Let $k\geq 5$.
By  (\ref{upper f_k}) and Theorem~\ref{k-critical}, for every $n\geq k,\, n\neq k+1$,
$$f_k(n+(k-1))-F(k,n+(k-1))\leq f_k(n)-F(k,n).$$
So the number of times when $f_k(n + k - 1) < f_k(n) + (k-1)(k - \frac{2}{k-1})/2 $ is bounded by
$$ \sum_{i = k+2}^{2k} f_k(n) - F(k,n). $$
Expanding (\ref{separation}), the above bound is at most
$$ \frac{1}{2}\sum_{i = k+2}^{2k - 2} \left(-i^2 + 3ik + \frac{k-3}{k-1}(k-i) - 2k^2 -2 \right) + 0 + \frac{k-2}{2} $$
$$ \leq \frac{-1}{12}\left(14k^3 - 45k^2 + 13k - 12\right) + \frac{9k^3 - 27k^2}{4} - \left(\frac{k^2 - 3k}{4} \cdot \frac{k-3}{k-1}\right) -k^3+3k^2 - k + 3 + \frac{k-2}{2} $$
$$ \leq \frac{k^3}{12} - \frac{k^2}{8} - \frac{11k}{6} + 7\leq  \frac{k^3}{12} - \frac{k^2}{8}. \qed $$

\section{Some applications}

\subsection{Ore-degrees}
\begin{figure}[htbp]
%\centering
\begin{center}
\includegraphics[height=4cm]{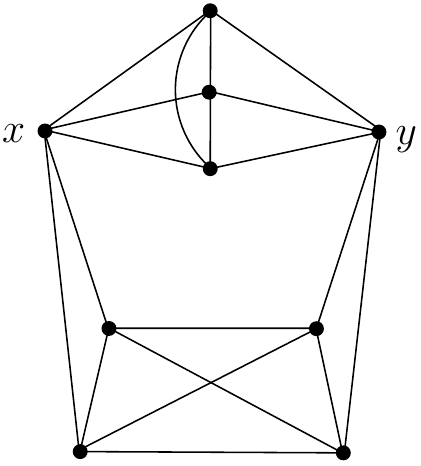}

\caption{The graph $O_5$.}
\label{fig:2}       % Give a unique label
\end{center}\end{figure}

The {\em Ore-degree}, $\Theta(G)$, of a graph $G$ is the maximum of $d(x)+d(y)$ over all edges $xy$ of $G$.
Let ${\mathcal{G}}_t=\{G\;:\; \Theta(G)\leq t\}$. It is easy to prove (see, e.g.~\cite{KK}) that $\chi(G)\leq 1+\lfloor t/2\rfloor$ for
every $G\in {\mathcal{G}}_t$. Clearly $\Theta(K_{d+1})=2d$  and $\chi(K_{d+1})=d+1$.
The graph $O_5$ in Fig~\ref{fig:2}  is the only $9$-vertex $5$-critical graph with $\Theta$ at most $9$.
We have $\Theta(O_5)=9$ and $\chi(O_5)=5$.

A natural question is to describe the graphs in ${\mathcal{G}}_{2d+1}$ with chromatic number $d+1$.
Kierstead and Kostochka~\cite{KK} proved that for $d\geq 6$ each such graph contains $K_{d+1}$.
Then Rabern~\cite{Ra} extended the result to $d=5$. Each $(d+1)$-chromatic graph $G$ contains a
$(d+1)$-critical subgraph $G'$. Since $\delta(G')\geq d$ and $\Theta(G')\leq \Theta(G)\leq 2d+1$,
\begin{equation}\label{kier}\mbox{
$\Delta(G')\leq d+1$, and  vertices of degree $d+1$ form an independent set. }
\end{equation}
Thus the results
in~\cite{KK} and~\cite{Ra} mentioned above could be stated in the following form.

\begin{theorem}[\cite{KK,Ra}]\label{d5} Let $d\geq 5$. Then the only $(d+1)$-critical graph $G'$
satisfying (\ref{kier}) is $K_{d+1}$.
\end{theorem}

The case $d=4$ was settled by Kostochka, Rabern, and Stiebitz~\cite{KRS}:

\begin{theorem}[\cite{KRS}]\label{d4} Let $d=4$. Then the only $5$-critical graphs $G'$
satisfying (\ref{kier}) are $K_{5}$ and $O_5$.
\end{theorem}

Theorem \ref{k-critical} and Corollary~\ref{co1} yield simpler proofs of Theorems~\ref{d5} and~\ref{d4}. The
key observation is the following.

\begin{lemma}\label{nl} Let $d\geq 4$ and $G'$ be a $(d+1)$-critical graph
satisfying (\ref{kier}). If $G'$ has $n$ vertices of which $h>0$ vertices have degree $d+1$, then
\begin{equation}\label{j1}
h\geq \left\lceil \frac{(d-2)n-(d+1)(d-2)}{d}\right\rceil
\end{equation}
and
\begin{equation}\label{j2}
h\leq \left\lfloor \frac{n-3}{d-1}\right\rfloor.
\end{equation}\end{lemma}

{\bf Proof.} By definition, $2e(G')=dn+h$. So, by Theorem \ref{k-critical} with $k=d+1$,
$$dn+h\geq (d+1 - \frac2{d})n- \frac{(d+1)(d-2)}{d},$$
which yields (\ref{j1}).

Let $B$ be the set of vertices of degree $d+1$ in $G'$ and $A=V(G')-B$.
By (\ref{kier}), $e(G'(A,B))=h(d+1)$. So, by Corollary~\ref{co1}(ii) with $k=d+1$,
$$h(d+1)\leq 3h+(n-h)-3=2h+n-3,$$
which yields (\ref{j2}).\qed

Another ingredient is the following old observation by Dirac.

\begin{lemma}[Dirac~\cite{D01}]\label{d1} Let $k\geq 3$. There are no $k$-critical graphs with $k+1$
vertices, and the only $k$-critical graph (call it $D_k$) with $k+2$
vertices is obtained from the $5$-cycle by adding $k-3$ all-adjacent vertices.
\end{lemma}

Suppose $G'$ with $n$ vertices of which $h$ vertices have degree
$d+1$ is a counter-example to Theorems~\ref{d5} or~\ref{d4}. Since the graph $D_{d+1}$
from Lemma~\ref{d1} has a vertex of degree $d+2$, $n\geq d+4$. So since $d\geq 4$, by~(\ref{j1}),
$$h\geq \left\lceil \frac{(d-2)(d+4)-(d+1)(d-2)}{d}\right\rceil=\left\lceil \frac{3(d-2)}{d}\right\rceil\geq 2.$$
On the other hand, if $n\leq 2d$, then by~(\ref{j2}),
$$h\leq \left\lfloor \frac{2d-3}{d-1}\right\rfloor=1.$$
Thus $n\geq 2d+1$.

Combining~(\ref{j1}) and~(\ref{j2}) together, we get
$$\frac{(d-2)n-(d+1)(d-2)}{d}\leq \frac{n-3}{d-1}.$$
Solving with respect to $n$, we obtain
\begin{equation}\label{j3}
n\leq \left\lfloor \frac{(d+1)(d-1)(d-2)-3d}{d^2-4d+2}\right\rfloor.
\end{equation}
For $d\geq 5$, the RHS of~(\ref{j3}) is less than $2d+1$, a contradiction to $n\geq 2d+1$.
This proves Theorem~\ref{d5}.

Suppose $d=4$. Then~(\ref{j3}) yields $n\leq 9$. So, in this case, $n=9$. By~(\ref{j1}) and~(\ref{j2}),
we get $h=2$. Let $B=\{b_1,b_2\}$ be the set of vertices of degree $5$ in $G'$.
By a theorem of Stiebitz~\cite{St},  $G'-B$ has at least two components. Since $|B|=2$ and $\delta(G')=4$,
each such component has at least $3$ vertices. Since $|V(G')-B|=7$, we may assume that
$G'-B$ has exactly two components, $C_1$ and $C_2$, and that $|V(C_1)|=3$. Again because $\delta(G')=4$,
$C_1=K_3$ and all vertices of $C_1$ are adjacent to both vertices in $B$. So, if we color both $b_1$ and
$b_2$ with the same color, this can extended to a $4$-coloring of $G'-V(C_2)$. Thus to have $G'$ $5$-chromatic,
we need $\chi(C_2)\geq 4$ which yields $C_2=K_4$. Since $\delta(G')=4$, $e(V(C_2),B)=4$. So, since each
of $b_1$ and $b_2$ has degree $5$ and $3$ neighbors in $C_1$, each of them has exactly two neighbors
in $C_2$. This proves Theorem~\ref{d4}.

\subsection{Local vs. global graph properties}

Krivelevich~\cite{Kr2} presented several nice applications of his lower bounds on $f_k(n)$ and related graph
parameters to questions of existence of complicated graphs whose small subgraphs are simple.
 We indicate here how to improve two of his bounds using Theorem~\ref{k-critical}.

Let $f(\sqrt{n},3,n)$ denote the maximum chromatic number over $n$-vertex graphs in which every
$\sqrt{n}$-vertex subgraph has  chromatic number at most $3$. Krivelevich proved that for every
fixed $\epsilon>0$ and sufficiently large $n$,
\begin{equation}\label{ju1}
f(\sqrt{n},3,n)\geq n^{6/31-\epsilon}.
\end{equation}
He used his result that every $4$-critical $t$-vertex graph with odd girth at least $7$ has
at least $31t/19$ edges. If instead of this result, we use  our bound on $f_4(n)$, then
repeating almost word by word Krivelevich's proof of his Theorem~4 (choosing $p=n^{-0.8-\epsilon'}$), we get  that for every
fixed $\epsilon$ and sufficiently large $n$,
\begin{equation}\label{ju2}
f(\sqrt{n},3,n)\geq n^{1/5-\epsilon}.
\end{equation}

\medskip
Another result of Krivelevich is:

\begin{theorem}[\cite{Kr2}]\label{kr3} There exists $C>0$ such that for every $s\geq 5$ there exists a
graph $G_s$ with at least $C\left(\frac{s}{\ln s}\right)^{\frac{33}{14}}$ vertices and  independence number less than $s$
such that the independence number of each $20$-vertex subgraph at least $5$.
\end{theorem}

He used the fact that for every $m\leq 20$ and every $m$-vertex $5$-critical graph $H$,
$$\frac{|E(H)|-1}{m-2}\geq \frac{\lceil 17m/8\rceil-1}{m-2}\geq\frac{33}{14}.$$
From  Theorem~\ref{k-critical} we instead get
$$\frac{|E(H)|-1}{m-2}\geq \frac{\left\lceil \frac{9m-5}{4}\right\rceil-1}{m-2}\geq\frac{43}{18}.$$
Then repeating the argument in~\cite{Kr2} we can replace $\frac{33}{14}$
in the statement of Theorem~\ref{kr3} with $\frac{43}{18}$.

\subsection{Coloring planar graphs}

One of the basic results on $3$-coloring of planar 
graphs is  Gr\" otzsch's Theorem~\cite{Gr}: every triangle-free planar graph is $3$-colorable.
The original proof of this theorem is somewhat sophisticated. 
There were subsequent simpler proofs (see, e.g.~\cite{Tho} and references
therein),
but Theorem~\ref{k-critical} yields a half-page proof.
A disadvantage of this  proof is that
 the proof of Theorem~\ref{k-critical} itself is not too simple.
In~\cite{KY}, we give a shorter proof of the fact $f_4(n) = F(4,n)$ and a short proof of Gr\" otzsch's Theorem.
In~\cite{BKY}, we use Theorem~\ref{k-critical} to 
 give  short proofs of some other known and new results on $3$-colorability of planar graphs.

\section{Algorithm}
Recall that $\rho_{k,G}(W)=(k+1)(k-2)|W|-2(k-1)|E(G[W])|$ and that $P_k(G)$ is the minimum of $\rho_{k,G}(W)$
over all nonempty $W\subseteq V(G)$. We will also use the related parameter $\widetilde{P}_k(G)$ which
is the minimum of $\rho_{k,G}(W)$
over all  $W\subset V(G)$ with $2\leq |W|\leq |V(G)|-1$.

\subsection{Procedure R1}
The input of the procedure $R1_k(G)$ is a  graph $G$. The output is one of the following five:\\
(S1) a nonempty set $R\subseteq V(G)$ with $\rho_{k,G}(R)\leq k(k-3)$, or\\
(S2) conclusion that $k(k-3)<\widetilde{P}_k(G)<(k+1)(k-2)$ and a nonempty  set $R\subsetneq V(G)$ with $\rho_{k,G}(R)= \widetilde{P}_k(G)$, or\\
(S3) conclusion that  $\widetilde{P}_k(G)<2(k-1)(k-2)$, and a set $R\subset V(G)$ with $2\leq |R|\leq n-1$
and  $\rho_{k,G}(R)=\widetilde{P}_k(G)$, or\\
(S4) conclusion that  $\widetilde{P}_k(G)=2(k-1)(k-2)$, and a set $R\subset V(G)$
 with $k\leq |R|\leq n-1$
and  $\rho_{k,G}(R)=2(k-1)(k-2)$, or\\
(S5) conclusion that  $\widetilde{P}_k(G)\geq 2(k-1)(k-2)$ and that every
set $R\subseteq V(G)$ with
$\rho_{k,G}(R)=2(k-1)(k-2)$ has size $k-1$ and induces $K_{k-1}$.

First we calculate $\rho_k(V(G))$, and if it is at most $k(k-3)$, then we are done. Suppose
\begin{equation}\label{9j}
(k+1)(k-2)|V(G)|-2(k-1)|E(G)|\geq 1+k(k-3).
\end{equation}

Consider the
auxiliary network $H=H(G)$ with vertex set $V\cup E\cup \{s,t\}$ and the set of
arcs $A=A_1\cup A_2\cup A_3$, where $A_1=\{sv\,:\,v\in V\}$,
$A_2=\{et\,:\,e\in E\}$, and $A_3=\{ve\,:\,v\in V,e\in E, v\in e\}$.
The capacity $c$ of each $sv\in A_1$ is $(k+1)(k-2)$, of each $et\in A_2$ is $2(k-1)$,
and of each $ve\in A_3$ is $\infty$.

Since the capacity of the cut $(\{s\}, V(H)-s)$ is finite, $H$ has
a maximum flow $f$.
Let $M(f)$ denote the value of  $f$, and
let $(S,T)$ be the minimum cut in it. By definition, $s\in S$ and $t\in T$.
Let $S_V=S\cap V$,  $S_E=S\cap E$, $T_V=T\cap V$, and $T_E=T\cap E$.

Since $c(ve)=\infty$ for every $v\in e$,
\begin{equation}\label{1}
\mbox{no edge of $H$ goes from $S_V$ to $T_E$.}
\end{equation}
It follows that if $e=vu$ in $G$ and $e\in T_E$, then
$v,u\in T_V$. On the other hand, if $e=vu$ in $G$,
$v,u\in T_V$ and $e\in S_E$, then moving $e$ from $S_E$
to $T_E$ would decrease the capacity of the cut by $2(k-1)$,
a contradiction. So, we get

\begin{claim}\label{cl1} $T_E=E(G[T_V])$.
\end{claim}

By the claim,
\begin{equation}\label{2}
M(f)=\min_{W\subseteq V}\bigg\{ (k+1)(k-2)|W|+2(k-1)(|E|-|E(G[W])|\bigg\}=2(k-1)|E|+\min\bigg\{P_k(G),0\bigg\}.
\end{equation}
So,  if $M(f)<2(k-1)|E|$, then $P_k(G)<0$ and any minimum cut gives us a set with small potential.
%it is not enough, since if $P_k(G)>0$, then the minimum cut will always have capacity $2(k-1)|E|$.
%But if $M(f)<2(k-1)|E|$, then we are done.
 Otherwise, consider for every $e_0\in E$ and
every vertex $v_0$ not incident to $e_0$, the network $H_{e_0,v_0}$ that has the same vertices and edges and
differs from
$H$ in the following:\\
(i) the capacity of the edge $e_0t$ is not $2(k-1)$ but $2(k-1)+2(k-1)(k-2)=2(k-1)^2$;\\
(ii) for every $v\in V(G)-v_0$, the capacity of the edge $sv$ is $(k+1)(k-2)-\frac{1}{2n}$;\\
(iii) the capacity of the edge $sv_0$ is $(k+1)(k-2)-\frac{1}{2n}+2(k-1)(k-2)+1$.\\

 Then for every $e_0\in E$ and $v_0\in V(G)$, the capacity
of the cut $(V(H_{e_0,v_0})-t,t)$ is $2(k-1)|E|+2(k-1)(k-2)$.
Since this is finite,   $H_{e_0,v_0}$ has a maximum flow  $f_{e_0,v_0}$.
As above, let $M(f_{e_0,v_0})$ denote the value of  $f_{e_0,v_0}$, and
let $(S,T)$ be the minimum cut in it. By definition, $s\in S$ and $t\in T$.
Let $S_V=S\cap V$,  $S_E=S\cap E$, $T_V=T\cap V$, and $T_E=T\cap E$.
By the same argument as above, (\ref{1}) and Claim~\ref{cl1} hold.
Let $M_k(G)$ denote the minimum value over $M(f_{e_0,v_0})$.

By (\ref{9j}), for every $e_0\in E$ and $v_0\in V(G)$, the capacity
of the cut $(s,V(H_{e_0,v_0})-s)$ is at least
$$\left((k+1)(k-2)-\frac{1}{2n}\right)n+2(k-1)(k-2)+1\geq 2(k-1)|E|+2(k-1)(k-2)+\frac{1}{2}.
$$

If the potential of some nonempty $W\neq V$ is less than $(k+1)(k-2)$, then
$G[W]$ contains some edge $e_0$ and there is $v_0\in V-W$. So,
in the network $H_{e_0,v_0}$, the capacity of the cut
$(\{s\}\cup (V-W)\cup (E-E(G[W])), W\cup E(G[W])\cup \{t\})$ is
$$\left((k+1)(k-2)-\frac{1}{2n}\right)|W|+2(k-1)(|E|-|E(G[W])|)=2(k-1)|E|+\rho_{k,G}(W)-\frac{|W|}{2n}.
$$
On the other hand, for every nonempty $W\neq V$,
every edge $e_0$ and every $v_0\in V$, the capacity of the cut
$(\{s\}\cup (V-W)\cup (E-E(G[W])), W\cup E(G[W])\cup \{t\})$ is at least
$$\left((k+1)(k-2)-\frac{1}{2n}\right)|W|+2(k-1)(|E|-|E(G[W])|)>2(k-1)|E|+\rho_{k,G}(W)-\frac{1}{2}.
$$
Thus if $M_k(G)\leq k(k-3)+2(k-1)|E|$, then (S1) holds and if $k(k-3)+2(k-1)|E|< M_k(G)<(k+1)(k-2)+2(k-1)|E|$,
 then (S2) holds.
Note that if a nonempty $W$ is independent, then $E(G[W])=\emptyset$, and the capacity of the cut
$(\{s\}\cup (V-W)\cup (E-E(G[W])), W\cup E(G[W])\cup \{t\})$ is at least
$$ 2(k-1)|E|+2(k-1)(k-2)+(k+1)(k-2).$$
Thus, if $$(k+1)(k-2)+2(k-1)|E| \leq M_k(G)<2(k-1)(k-2)-1+2(k-1)|E|,$$
 then (S3) holds.

Similarly, if
 $$2(k-1)(k-2)-1+2(k-1)|E| \leq M_k(G)<2(k-1)(k-2)+2(k-1)|E|-\frac{k-1}{2n},$$
 then there exists  $W\subset V$ with $k\leq |W|\leq n-1$
 with potential $2(k-1)(k-2)$. Then (S4) holds. Finally, if
$M_k(G)\geq 2(k-1)(k-2)+2(k-1)|E|-\frac{k-1}{2n},$
then (S5) holds.

Since the complexity of the max-flow problem is at most $Cn^2\sqrt{|E|}$ and $|E|\leq kn$, the procedure takes time
at most $Ck^{1.5}n^{4.5}$.

\subsection{Outline of the algorithm}
We consider the outline for $k\geq 7$. For $k\leq 6$, everything is quite similar and easier.

 Let the input be an $n$-vertex $e$-edge graph $G$.
  The algorithm will be recursive.
 The output will be either a coloring of $G$ with $k-1$ colors or return a nonempty $R\subseteq V(G)$ with $\rho_{k,G}(R)\leq k(k-3)$.
The algorithm runs through $7$ steps, which are listed below.
If a step is triggered, then a recursive call is made on a smaller graph $G'$.
Some steps will then require a second recursive call on another graph $G''$.

The algorithm does not make the recursive call if $|E(G')| \leq k^2 / 2$.
%If $G'$ is isomorphic to $K_{k-1}$, then the coloring is made directly.
%If $G'$ is not isomorphic to $K_{k-1}$, then $G'$ is $(k-2)$-degenerate.
In this case, $G'$ is either $(k-2)$-degenerate or $K_k$ minus a matching, and so is easily $(k-1)$-colorable in time $O(k|V(G')|^2)$.
%the algorithm colors $G'$ using standard methods for degenerate graphs.
This also holds for $G''$.

After all calls have been made, the algorithm will return a coloring or a subgraph with low potential, skipping the other steps.

1) We check whether $G$ is disconnected or has a cut-vertex or has a vertex of degree at most~$k-2$.
In the case of any "yes", we consider smaller graphs (and at the end will reconstruct the coloring).\\

2) We run  $R1_k(G)$ and consider possible outcomes. If the outcome is (S1), we are done.\\

3) Suppose the outcome is (S2).
The algorithm makes a recursive call on $G' = G[R]$, which returns a $(k-1)$-coloring $\phi$.
Let $G''$ be the graph $Y(G,R,\phi)$ described in Definition~\ref{def1}.
The proof of Claim~\ref{very small potential} yields that $P_k(G'')\geq k(k-3)$, and thus the recursive call will return with a coloring.
Let $\phi'$ be the coloring returned.
It is straightforward to combine the colorings $\phi$ and $\phi'$ into a $(k-1)$-coloring of $G$.\\

4) Suppose the outcome is (S3) or (S4). We choose $i$ using (\ref{j14}) and add $i$
edges to $G[R]$ as in the proof of Claim~\ref{small potential}. Denote the new graph $G'$.
 The algorithm makes a recursive call on $G' = G[R]$, which returns a $(k-1)$-coloring $\phi$.
  Let $G''$ be the graph $Y(G,R,\phi)$ described in Definition~\ref{def1}. The proof of
Claim~\ref{small potential} yields that $P_k(G'')\geq k(k-3)$, and thus the recursive call will return with a coloring.
Let $\phi'$ be the coloring returned.
It is straightforward to combine the colorings $\phi$ and $\phi'$ into a $(k-1)$-coloring of $G$.\\

5) So, the only remaining possibility is (S5). For every $(k-1)$-vertex $v\in V(G)$, check whether
there is a $(k-1)$-clique $K(v)$ containing $v$ (since (S5) holds, such a clique is unique, if exists).
 We certainly can do this in $O(kn^2)$ time.
 Let $a_v$ denote the neighbor of $v$ not in $K(v)$ and $T_v$ denote the set of $(k-1)$-vertices in $K(v)$.
Then for every pair $(v,K(v))$ such that $d(v)=k-1$ and $K(v)$ exists, do the following:\\

(5.1) If there is  $w\in T_v- v$ with $a_w\neq a_v$, then consider
the graph $G'=G-v-w+a_va_w$. By Claim~\ref{clusters in cliques}, $P_k(G')>k(k-3)$.
So, the algorithm will return with a $(k-1)$-coloring of $G'$, which we then extend to $G$.

(5.2) Suppose that $|T_v|\geq 2$ and $K(v)-T_v$ contains a vertex $x$ of degree at most $k-2+|T_v|$.
Let $G'=G-x+v'$, where the closed neighborhood of $v'$ is the same as of $v$.
By Claim~\ref{big neighbors}, $P_k(G')>k(k-3)$, so the algorithm returns a $(k-1)$-coloring of $G'$, which is then extended to $G$ as in the proof of Claim~\ref{big neighbors}.

(5.3) Suppose that $T_v=\{v\}$ and $K(v)$ contains at least $k/2 - 1$ vertices
of degree $k$. Since (S5) holds, there is $x\in K(v)-v$ of degree at most $k$ not adjacent to $a_v$.
Let $x_1$ and $x_2$ be the neighbors of $x$ outside of $K_v$.
Let $G'$ be obtained from $G-v$ by adding edges $a_vx_1$ and $a_vx_2$.
By the proof of Claim~\ref{cliques have k+1},  $P_k(G')>k(k-3)$, so the algorithm finds a $(k-1)$-coloring of $G'$, which is then extended to $G$ as in the proof of Claim~\ref{cliques have k+1}.\\

6) Let $C_v$ denote the {\em cluster of $v$}, i.e.  the set of vertices that have the same closed neighborhood as $v$.
We certainly can find $C_v$ for every  $(k-1)$-vertex $v\in V(G)$  in $O(kn^2)$ time.
Then for every pair $(v,C_v)$
 such that $d(v)=k-1$, do the following:\\

(6.1) Suppose that $|C_v| \geq 2$ and $N(v)-C_v$ contains a vertex $x$ of degree at most $k-2+|T_v|$.
Then do the same as in (5.2).

(6.2) Suppose that $N(v)-C_v$ contains a  $(k-1)$-vertex $w$ and that $|C_w|\leq |C_v|$.
If $v$ is not in a $(k-1)$-clique, then consider $G'=G-w+v'$, where the $v'$ is a new vertex
whose closed neighborhood is the same as that of $v$.
By the proof of Claim~\ref{adjacent k-1}, $P_k(G')>k(k-3)$, and so we find a $(k-1)$-coloring of $G'$ and then extend it to $G$
as in the proof of Claim~\ref{adjacent k-1}.\\

7) Let $L_0$, $H_0$, and $e_0$ be as defined in Definition \ref{L,H,e0}.
If $e_0 \geq 2(|L_0| + |H_0|)$, then iteratively remove vertices in $L_0$ with at most two neighbors in $H_0$ and vertices in $H_0$ with at most two neighbors in $L_0$.
Let $H$ be the graph that remains, and $G' = G - V(H)$.
Clearly $P_k(G') > k(k-3)$, so the recursive call returns a coloring of $G'$.
Give each vertex $v \in V(H)$ a list of colors  $L(v) = \{c_1, \dots, c_{k-1}\}$, then remove from that list the colors on $N(v) \cap V(G')$.
Orient the edges of $H$ as in Case $1$ of the proof of Lemma \ref{lem3}.
Then extend the coloring of $G'$ to a coloring of $G$ by list coloring $H$ using
the system described in the proof to Lemma~\ref{lem1}.

\subsection{Analysis of correctness and running time}
The proof of Theorem \ref{k-critical} consists in proving that at least one of the situations in steps $1$ through $7$ described above must happen.
Moreover, the main theorem proves that $G', G'' \prec G$ by a partial order with finite descending chains, and therefore the algorithm will terminate.
We claim that the algorithm makes at most $O(k^2 n^2  \log(n))$ recursive calls, and each call only takes $O(k^{1.5}n^{4.5})$ time, so the algorithm runs in $O(k^{3.5}n^{6.5}\log(n))$ time.

If a call of the recursive algorithm terminates on Step $2$, we will refer to this as `Type 1,' a call terminating on Step $1$, $3$, $4$, $5.1$, $5.3$, $6.1$, or $7$ is `Type 2,' and a call terminating on Step $5.2$ or $6.2$ is `Type 3.'
If a call is made on a Type 1, then the whole algorithm stops.

If a Type 3 happens, then the algorithm makes one recursive call with a graph with the same number of edges and strictly more pairs of vertices with the same closed neighborhood.
The proof of Claim \ref{small clusters} shows that the number of pairs of vertices with the same closed neighborhood is bounded by $kn$.
Then at least one out of every $kn$ consecutive recursive calls is Type 1 or 2.

Consider an instance of a Type 2 call with input graph $H$.
If $H'$ is the graph in the first recursive call and $H''$ is the graph in the second call (if necessary), then $|E(H')|, |E(H'')| < |E(H)|$ and $|E(H)| \geq |E(H')| + |E(H'')| - k^2 / 2$.
Let $g_k(e, i)$ denote the number of Type 2 recursive calls made on graphs with $i$ edges.
Note that if $i \leq k^2/2$ then $g_k(e, i) = 0$ and $g_k(e,e) = 1$.
By tracing calls up through their parent calls, it follows that
$$	e \geq i + (g_k(e, i) - 1)\left(i - k^2/2 \right)
$$
when $i > k^2 / 2$.
Therefore
$$ g_k(e, i) < \frac{e}{ \left(i - k^2/2  \right)}.
$$
The total number of calls that our algorithm makes is at most
$$ k n \sum_{i = k^2/2 + 1} ^e g_k(e, i) < k n e \log(e).
$$
Because $e \leq nk$, we have that the total number of calls is $O(k^2 n^2  \log(n))$.

A call may run algorithm R1 once, which will take $O(k^{1.5}n^{4.5})$ time.
Constructing the appropriate graphs for recursion in Steps $3$, $4$, $5$, and $6$ will take $O(kn^2)$ time.
Combining colorings in Steps $1$, $3$, $4$, $5$, and $6$ will take $O(n)$ time.
Coloring a degenerate graph will take $O(k n^2)$ time, which happens at most twice.
The only thing left to consider is Step $7$.
Iteratively removing vertices will take $O(n^2)$ time.
Splitting the vertices and orienting the edges using network flows will take $O(n^{2.5}k^{0.5})$ time.
Finding a kernel will take $O(n^2)$ time, which happens at most $n$ times.
Therefore each instance of the algorithm takes $O(k^{1.5}n^{4.5})$ time.

\section{Concluding remarks}

Many open questions remain.

1. It would be good to find exact values of $f_k(n)$ for all $k$ and $n$.

2. Similar questions for list coloring look much harder. Some results
are in~\cite{K3,K5}.

3. One can ask how few edges may have an $n$-vertex $k$-critical graph
not containing a given subgraph, for example, with bounded clique number.
Krivelevich~\cite{Kr2} has interesting results on the topic.

4. Brooks-type results would be interesting.

5. A similar problem for hypergraphs was considered in~\cite{K3}, but the bounds there are good
only for large $k$.

6. It is clear that there are algorithms with better performance than ours.

\bigskip
{\bf Acknowledgment.} The authors thank Xuding Zhu for the nice reduction idea and Oleg Borodin,
Michael Krivelevich, Bernard Lidicky
and Artem Pyatkin for helpful comments.


\begin{thebibliography}{99}

\bibitem{Bor}
O. V. Borodin,  A new proof of Gr\" unbaum's 3 color theorem, {\em Discrete Math.}
{\bf 169} (1997),  177--183.

\bibitem{BKY}
O. V. Borodin, A.~V.~Kostochka, B. Lidicky, and M. Yancey, Short proofs of coloring theorems on
planar graphs, in preparation.

\bibitem{Br} R.~L.~Brooks, On colouring the nodes of a network,
    {\em Math. Proc.
    Cambridge Philos. Soc.} {\bf 37} (1941), 194--197.

\bibitem{D0}
 G.~A.~Dirac, Note on the colouring of graphs, {\it Math. Z.}
 {\bf 54} (1951), 347-353.

\bibitem{D02}
 G.~A.~Dirac,
A property of 4-chromatic graphs and some remarks on critical graphs,
 J. London Math. Soc. 27 (1952),
 85-92.

\bibitem{D03}
 G.~A.~Dirac,
 Some theorems on abstract graphs, Proc. London Math. Soc. (3) 2 (1952), 69-81.

\bibitem{D01}
 G.~A.~Dirac, Map colour theorems related to the Heawood colour formula, {\it J. London
 Math. Soc.} {\bf 31} (1956), 460-471.

\bibitem{D1}
 G.~A.~Dirac, A theorem of R.~L.~Brooks and a conjecture of H.~Hadwiger,
 {\it Proc. London Math. Soc.} (3) {\bf 7} (1957), 161-195.

\bibitem{D2}
 G.~A.~Dirac, The number of edges in critical graphs, {\it J. Reine Angew.
 Math.} {\bf 268/269} (1974), 150-164.


\bibitem{FM}
B. Farzad, M. Molloy,  On the edge-density of 4-critical graphs,
{\em Combinatorica} {\bf 29} (2009),  665--689.


\bibitem{G1}
 T.~Gallai, Kritische Graphen I, {\it Publ. Math. Inst. Hungar. Acad. Sci.}
 {\bf  8} (1963), 165-192.

\bibitem{G2}
 T.~Gallai, Kritische Graphen II, {\it Publ. Math. Inst. Hungar. Acad. Sci.}
 {\bf 8} (1963), 373-395.

\bibitem{Gr}
H. Gr\" otzsch,  Zur Theorie der diskreten Gebilde. VII. Ein Dreifarbensatz f\" ur dreikreisfreie Netze auf der Kugel.
{\em Wiss. Z. Martin-Luther-Univ. Halle-Wittenberg. Math.-Nat. Reihe} {\bf 8} (1958/1959), 109--120 (in German).

\bibitem{H}
 G.~Haj\'os, \"Uber eine Konstruktion nicht $n$--f\"arbbarer Graphen,
 {\it Wiss. Z. Martin-Luther-Univ. Halle-Wittenberg Math.-Naturw. Reihe}
 {\bf  10} (1961), 116--117.


\bibitem{J}
 T.~R.~Jensen and B.~Toft, Graph Coloring Problems, Wiley-Interscience
 Series in Discrete Mathematics and Optimization, John Wiley \& Sons,
 New York, 1995.

\bibitem{J2}
 T.~R.~Jensen and B.~Toft,
 25 pretty graph colouring problems, {\em Discrete Math.} {\bf 229} (2001),  167-169.

\bibitem{KK} H.~A.~Kierstead, A.~V.~Kostochka, Ore-type versions of
    Brook's theorem, {\em Journal of Combinatorial Theory, Series
    B.} {\bf 99} (2009), 298--305.

\bibitem{KRS}
 A.~V.~Kostochka, L. Rabern, and M.~Stiebitz, Graphs with chromatic number close to maximum degree,
 to appear in {\em Discrete Math.}

\bibitem{K2}
 A.~V.~Kostochka and M.~Stiebitz, Excess in colour-critical graphs, in:
Graph Theory and Combinatorial Biology, Balatonlelle (Hungary), 1996,
{\it Bolyai Society, Mathematical Studies} {\bf  7}, Budapest, 1999, 87-99.

\bibitem{K3}
 A.~V.~Kostochka and M.~Stiebitz, On the number of edges in colour-critical
graphs and hypergraphs, {\em Combinatorica} {\bf 20} (2000), 521--530.


\bibitem{K4}
A.~V.~Kostochka and M.~Stiebitz, A list version of Dirac's theorem on the
number of edges in colour-critical graphs, {\em Journal of Graph Theory}
{\bf 39} (2002), 165--167.


\bibitem{K5}
A.~V.~Kostochka and M.~Stiebitz, A new lower bound on the number of edges
in colour-critical graphs and hypergraphs, {\em Journal of Combinatorial Theory, Series
B.} {\bf 87} (2003), 374--402.

\bibitem{KY}
A.~V.~Kostochka and M. Yancey, Ore's Conjecture for $k=4$ and Gr\" otzsch Theorem, in preparation.

\bibitem{Kr2}
 M.~Krivelevich, On the minimal number of edges in color-critical graphs,
 {\it Combinatorica} {\bf 17} (1997), 401-426.

%\bibitem{Kro}
% H.~V.~Kronk and J.~Mitchem, On Dirac's generalization
% of Brooks' theorem, {\it Canad. J. Math.} {\bf 24} (1972), 805-807.

%\bibitem{L}
% L.~Lov\'asz, Combinatorial Problems and Exercises, Akad. Kiad\'o, Budapest,
% 1979.

%\bibitem{M}
% J.~Mitchem, A new proof of a theorem of Dirac on the
% number of edges in critical graphs,  {\it J. Reine u. Angew.
% Math.} {\bf 299/300} (1978), 84-91.

\bibitem{O}
 O.~Ore, The Four Color Problem, Academic Press, New York, 1967.


\bibitem{Ra} L.~Rabern, Coloring $\Delta$-critical graphs with small
    high vertex cliques, manuscript 2010.

\bibitem{Ra1} L.~Rabern, An improvement on Brooks's theorem,
    manuscript 2010.

\bibitem{St} M.~Stiebitz, Proof of a conjecture of T. Gallai
    concerning connectivity properties of
    colour-critical graphs, {\em Combinatorica} {\bf 2} (1982), 315--323.

\bibitem{Tho}
C. Thomassen,  A short list color proof of Gr\" otzsch's theorem,
{\em J. Combin. Theory Ser. B} {\bf 88} (2003),  189--192.


%\bibitem{Toft1}
% B.~Toft,
%75 graph-colouring problems. Graph colourings (Milton Keynes, 1988), 9--35, Pitman Res. Notes Math. Ser.,
% 218, Longman Sci. Tech., Harlow, 1990.

\bibitem{Tuza1} Zs. Tuza, Graph coloring, in:
Handbook of graph theory.
 J. L. Gross and J. Yellen Eds.,  CRC Press, Boca Raton, FL, 2004. xiv+1167 pp. 

%\bibitem{West}
%D. B. West, Introduction to Graph Theory.  Second edition Prentice Hall, 2001.

\end{thebibliography}
\end{document}